\documentclass{article}
\usepackage{amsmath,amssymb,indentfirst}
\usepackage[latin1]{inputenc}
\usepackage{url}
\newtheorem{theorem}{Theorem}
\newtheorem{lemma}{Lemma}

\begin{document}

\title{THE IMPLICIT AND THE INVERSE FUNCTION THEOREMS: EASY PROOFS}

\author{Oswaldo Rio Branco de Oliveira}
\date{}
\maketitle

\begin{abstract} This article presents simple and easy proofs of the Implicit Function Theorem and the Inverse Function Theorem, in this order, both of them on a finite-dimensional Euclidean space, that employ only the Intermediate Value Theorem and the Mean-Value Theorem. These proofs avoid compactness arguments, the contraction principle, and fixed-point theorems.
\end{abstract}

\vspace{0,2 cm}

\hspace{- 0,6 cm}{\sl Mathematics Subject Classification: 26B10, 26B12, 97I40, 97I60}

\hspace{- 0,6 cm}{\sl Key words and phrases:} Implicit Function Theorems, Calculus of Vector Functions, Differential Calculus, Functions of Several Variables.

\vspace{0,3 cm}

\section{Introduction.} The objective of this paper is to present very simple and easy proofs of the Implicit and Inverse Function theorems, in this order, on a finite-dimensional Euclidean space. Following Dini's approach, these proofs do not employ compactness arguments, the contraction principle or any fixed-point theorem. Instead of such tools, the proofs presented in this paper rely on the Intermediate Value Theorem and the Mean-Value Theorem on the real line.

The history of the Implicit Function Theorem is quite long and dates back to R. Descartes (on algebraic geometry), I. Newton, G. Leibniz, J. Bernoulli, and L. Euler (and their works on infinitesimal analysis), J. L. Lagrange (on real power series),  A. L. Cauchy (on complex power series), to whom is attributed the first rigorous proof of such result, and U. Dini (on functions of real variables and differential geometry). For extensive accounts on the history of the Implicit Function Theorem (and the Inverse Function Theorem) and further developments (as in differentiable manifolds, Riemannian geometry, partial differential equations, numerical analysis, etc.) we refer the reader to Krantz and Parks \cite{Krantz}, Dontchev and Rockafellar [7, pp. 7--8, 57--59], Hurwicz and Richter \cite{Hurwicz}, 
Iurato \cite{Iurato}, and Scarpello \cite{Scarpello}. 

Complex variables versions of the theorems studied in this paper can be seen in Burckel [4, pp. 173--174, 180--183] and in Krantz and Parks [14, pp. 27--32, 117--121]. See also de Oliveira [5, p. 17].

The two most usual approaches to the Implicit and Inverse Function theorems on a finite-dimensional Euclidean space begin with a proof of the Inverse Function Theorem and then derive from it the Implicit Function Theorem. The most basic of these approaches employs the compactness of the bounded and closed balls inside $\mathbb R^n$, where $n\geq 1$, and Weierstrass' Theorem on Minima. See Apostol [1, pp. 139-148], Bartle [2, pp. 380--386], Fitzpatrick [8, pp. 433--453], Knapp [13, pp. 152--161], Rudin [17, pp. 221--228], and Spivak [19, pp. 34--42]. The more advanced approach is also convenient on complete normed vector spaces (also called Banach spaces) of infinite dimension, employs the contraction mapping principle (also called the fixed-point theorem of Banach), and besides others applications is very useful in studying some types of partial differential equations. See Dontchev and Rockafellar [7, pp. 10--20], Gilbarg and Trudinger [9, pp. 446--448], Renardy and Rogers [16, pp. 336--339], and Taylor [20, pp. 11--13]. 


In this article, we prove at first the Implicit Function Theorem, by induction, and then we derive from it the Inverse Function Theorem. This approach is accredited to U. Dini (1876), who was the first to present a proof (by induction) of the Implicit Function Theorem for a system with several equations and several real variables, and then stated and also proved the Inverse Function Theorem. See Dini [6, pp. 197--241].


Another proof by induction of the Implicit Function Theorem, that also simplifies Dini's argument, can be seen in the book by Krantz and Parks [14, pp. 36--41].
Nevertheless,  this particular proof by Krantz and Parks does not establish the local uniqueness of the implicit solution of the given equation. On the other hand, the proof that we present in this paper further simplifies Dini's argument and  makes the whole proof of the Implicit Function Theorem  very simple, easy, and with a very small amount of computations. The Inverse Function Theorem then follows immediately.

\section{Notations and Preliminaries.}

\vspace{0,2 cm}

Let us indicate by $\mathbb R$ the complete field of the real numbers, and by $\mathbb R^n$, where $n\geq 1$, the $n$-dimensional Euclidean space. If $a$ and $b$ are real numbers, with $a<b$, we set $]a,b[=\{x\ \textrm{in}\ \mathbb R: a<x<b\}$ and $[a,b]=\{x \ \textrm{in}\ \mathbb R: a\leq x\leq b\}$. We assume without proof the following two results.
\begin{itemize}
\item[$\bullet$] {\sf The Intermediate Value Theorem:} Let $f:[a,b]\to \mathbb R$ be continuous. If $\lambda$ is a value between $f(a)$ and $f(b)$, then there is a $c$ in $[a,b]$ satisfying $f(c)=\lambda$.
\item[$\bullet$] {\sf The Mean-Value Theorem:} Let $f:[a,b]\to \mathbb R$ be continuous on $[a,b]$ and differentiable on the open interval $]a,b[$. Then, there exists $c$ in $]a,b[$ satisfying $f(b)-f(a)=f'(c)(b-a)$.
\end{itemize}
The first complete proof of the Intermediate Value Theorem was given by B. Bolzano (1817), who was also the first to enunciate the supremum property of the real numbers, see Lützen [15, p. 175] and Hairer and Wanner [10, p. 206]. The Mean-Value Theorem is attributed to Lagrange (1797), see Hairer and Wanner [10, p. 240].

In what follows we fix the ordered canonical bases $\{e_1,\ldots ,e_n\}$ and $\{f_1, \ldots ,f_m\}$, of $\mathbb R^n$ and $\mathbb R^m$, respectively. Given $x=x_1e_1 + \cdots + x_ne_n$ in $\mathbb R^n$, where $x_1,\ldots,x_n$ are real numbers, we write $x=(x_1,\ldots,x_n)$.
If $x=(x_1,\ldots, x_n)$ and $y=(y_1,\ldots,y_n)$ are in $\mathbb R^n$, then their inner product is defined by $\left<x, y\right>= x_1y_1 + \cdots + x_ny_n$. The norm of $x$ is $|x|=\sqrt{x_1^2 + \cdots +x_n^2}=\sqrt{\left<x,x\right>}$. It is well-known the Cauchy-Schwarz inequality $|\left<x,y\right>|\leq |x||y|$, for all $x$ and all $y$, both in $\mathbb R^n$. Moreover, given a point $x$ in $\mathbb R^n$ and $r>0$, the set $B(x;r)=\{y\ \textrm{in}\ \mathbb R^n:\ |y-x|<r\}$ is the open ball centered at $x$ with radius $r$.

Given a linear function $T:\mathbb R^n \rightarrow \mathbb R^m$, there are real numbers $a_{ij}$, where $i=1,\ldots,m$ and $j=1,\ldots,n$, such that $T(e_j)= a_{1j}f_1 + \cdots + a_{mj}f_m$, for each $j=1,\ldots,n$. We associate to $T$ the $m\times n$ real matrix 
\[\left(\begin{array}{lllll}
a_{11} & \ldots &   a_{1n} \\
\ \vdots &      & \ \vdots \\
a_{m1} & \ldots &   a_{mn}\\
\end{array}
\right)\ = \ (a_{ij})_{\substack{1\leq i\leq m\\ 1\leq j\leq n}}.\]
The jth column of the matrix $(a_{ij})$ supplies the $m$ coordinates of $T(e_j)$. Furthermore, given a $m\times n$ real matrix $M=(a_{ij})$ we associate with it the linear function $T: \mathbb R^n\to \mathbb R^m$ defined by $T(e_j)= (a_{1j},\ldots,a_{mj})$, for all $j=1,\ldots,m$. If $v$ belongs to $\mathbb R^n$, for simplicity of notation we also write $Tv$ for $T(v)$. The Hilbert-Schmidt norm of $T$ is the number $\|T\|=\sqrt{\sum_{i=1}^m\sum_{j=1}^n|a_{ij}|^2}$.

\begin{lemma}\label{L1} Let $T:\mathbb R^n \rightarrow \mathbb R^m$ be a linear function. We have $|Tv|\leq \|T\||v|$, for all $v$ in $\mathbb R^n$, and $T$ is continuous. 

\end{lemma}

\hspace{-0,6 cm}{\bf Proof.} Let $(a_{ij})$ be the matrix associated to $T$. Given $v$ in $\mathbb R^n$, we have $|Tv|^2= \sum_{i=1}^m|\left<(a_{i1},\ldots,a_{in}),v\right>|^2$. By employing the Cauchy-Schwarz inequality we obtain $|Tv|^2 \leq \sum_{i=1}^m |(a_{i1},\ldots,a_{in})|^2|v|^2= |v|^2\sum_{i=1}^m\sum_{j=1}^n|a_{ij}|^2=\|T\|^2|v|^2$. Thus, $|Tv|\leq \|T\||v|$. Finally, from the inequality $|T(v+h)- T(v)|=|Th|\leq \|T\||h|$, for all $v$ and all $h$, both in $\mathbb R^n$, it follows that $T$ is continuous on $\mathbb R^n\blacksquare$


\vspace{0,2 cm}

Let $\Omega$ be an open set in $\mathbb R^n$. Given a function $F: \Omega \rightarrow \mathbb R^m$ and $p$ in $\Omega$, we write $F(p)=\big(F_1(p),\ldots,F_m(p)\big)$, with  $F_i:\Omega\rightarrow \mathbb R$ the ith component of $F$, for each $i=1,\ldots, m$. We say that $F$ is differentiable at $p$ if there is a linear function $T:\mathbb R^n \rightarrow R^m$ and a function $E:B(0;r)\rightarrow \mathbb R^m$ defined on some ball $B(0;r)$, with $r>0$, centered at the origin of $\mathbb R^n$ such that  
$F(p+h) = F(p) + T(h) + E(h)$, for all $h$ in $B(0;r)$, where $E(h)/|h| \to 0$ as $h\to 0$.
The function $F$ is differentiable if it is differentiable at every point in $\Omega$. 

The linear map $T$ is the differential of $F$ at $p$, denoted by $DF(p)$. Supposing that $t$ is a real variable and defining for each $j=1,\ldots, n$ the directional derivative of $F$ at $p$ in the direction $e_j$ as $\frac{\partial F}{\partial e_j}(p)=\lim\limits_{t\to 0}\frac{F(p+te_j) - F(p)}{t}$, we deduce that $DF(p)(e_j)=\frac{\partial F}{\partial e_j}(p)=\big(\frac{\partial F_1}{\partial e_j}(p),\ldots,\frac{\partial F_m}{\partial e_j}(p)\big)$. Putting $\frac{\partial F_i}{\partial x_j}(p)=\frac{\partial F_i}{\partial e_j}(p)$, for all $i=1,\ldots, m$ and all $j=1,\ldots, n$, we call $\frac{\partial F}{\partial x_j}(p)=\big(\frac{\partial F_1}{\partial x_j}(p),\ldots,\frac{\partial F_m}{\partial x_j}(p)\big)$ the jth partial derivative of $F$ at $p$. The Jacobian matrix of $F$ at $p$ is 
\[JF(p)=\left(\frac{\partial F_i}{\partial x_j}(p)\right)_{\substack{1\leq i\leq m\\1 \leq j\leq n}}=
\left(\begin{array}{lllll}
\frac{\partial F_1}{\partial x_1}(p) & \cdots & \frac{\partial F_1}{\partial x_n}(p)\\
\ \ \ \vdots  &   & \ \ \ \vdots \ \ \      \\
\frac{\partial F_m}{\partial x_1}(p) & \cdots & \frac{\partial F_m}{\partial x_n}(p)
\end{array}
\right).\]
We say that $F$ is of class $C^1$ if $F$ and its partial derivatives of order $1$ are continuous on $\Omega$. In such case, we also say that $F$ is in $C^1(\Omega;\mathbb R^n)$.

The following lemma is a local result and valid on an open subset $\Omega$ of $\mathbb R^n$. For practicality, we enunciate it for $\Omega=\mathbb R^n$. 

\begin{lemma}\label{L2} Let $F:\mathbb R^n \to \mathbb R^m$ be  differentiable, $T: \mathbb R^k \to \mathbb R^n$ be the linear function associated to a $n\times k$ real matrix $M$, and $y$ be a fixed point in $\mathbb R^n$. Then, the function $G(x)= F(y + Tx)$, where $x$ is in $\mathbb R^k$, is differentiable and satisfies $JG(x)=JF(y + Tx)M$, for all $x$ in $\mathbb R^k$.
\end{lemma}

\hspace{-0,6 cm}{\bf Proof.} Let us fix $x$ in $\mathbb R^k$. Given $v$ in $\mathbb R^n$, by the differentiability of $F$ we have
$F(y+ Tx + v) = F(y+ Tx) + DF(y+ Tx)v + E(v)$, where $E(v)/|v|\rightarrow 0$ as $v \rightarrow 0$.
Substituting $v=Th$, where $h$ is in $\mathbb R^k$, into the last identity we obtain $ G(x+h)= G(x) + DF(y+ Tx)Th + E(Th)$. Thus, supposing $h\neq 0$, we have $\frac{E(Th)}{|h|}= 0$, if $Th=0$, and $\frac{E(Th)}{|h|}=\frac{E(Th)}{|Th|}\frac{|Th|}{|h|}$, if $Th\neq 0$. Since $Th\rightarrow 0$ as $h\rightarrow 0$ and, by Lemma 1, $|Th|/|h|\leq \|T\|$, we conclude that $E(Th)/|h|$ goes to $0$ as $h$ goes to $0$. Hence, $G$ is differentiable at $x$ and $JG(x)= JF(y+Tx)M\blacksquare$

\vspace{0,10 cm}

With the hypothesis on Lemma \ref{L2}, we see that if $F$ is $C^1$ then $G$ also does. Given $a$ and $b$, both in $\mathbb R^n$, we put  $\overline{ab}=\{a+t(b-a):0\leq t\leq 1\}$.


\begin{lemma}\label{L3}{\bf (The Mean-Value Theorem in Several Variables)}. Let us consider $F:\Omega \to \mathbb R$ differentiable, with $\Omega$ open in $\mathbb R^n$. Let $a$ and $b$ be points in $\Omega$ such that the segment $\overline{ab}$ is within $\Omega$. Then, there exists $c$ in $\overline{ab}$ satisfying
\[F(b)-F(a) = \left<\nabla F(c), b-a\right>.\]
\end{lemma}
{\bf Proof.} The curve $\gamma(t)= a +t(b-a)$, with $t$ in $[0,1]$, is inside $\Omega$. Employing the mean-value theorem in one real variable and Lemma \ref{L2} we obtain
$F(b)-F(a)= (F\circ\gamma)(1)-(F\circ\gamma)(0)=(F\circ\gamma)'(t_0)=\left<\nabla F(\gamma(t_0)), b-a \right>$, for some $t_0$ in $[0,1]\blacksquare$


\begin{lemma} \label{L4}Let $F$ be in $C^1(\Omega; \mathbb R^n)$, with $\Omega$ open within $\mathbb R^n$, and $p$ in $\Omega$ satisfying $\det JF(p)\neq 0$. Then, $F$ restricted to some ball $B(p;r)$, with $r>0$, is injective.
\end{lemma}

\hspace{- 0,6 cm}{\bf Proof.} 
(See Bliss \cite{Bliss})
 Since $F$ is of class $C^1$ and the determinant function $\det: \mathbb R^{n^2}\to \mathbb R$ is continuous and $\det JF(p)= \det \big(\frac{\partial F_i}{\partial x_j}(p)\big)\neq 0$, there is $r>0$ such that $\det\big(\frac{\partial F_i}{\partial x_j}(\xi_{ij})\big)$ does not vanish, for all $\xi_{ij}$ in $B(p;r)$, where $1\leq i,j\leq n$. 

Let $a$ and $b$ be distinct in $B(p;r)$. By employing 
the mean-value theorem in several variables 
to each component $F_i$ of $F$, we find $c_i$ in the segment $\overline{ab}$, within  $B(p;r)$, such that $F_i(b)-F_i(a)=\left<\nabla F_i(c_i),b-a\right>$. Hence, 
\[ \left(\begin{array}{l}
F_1(b)- F_1(a)\\
\ \ \ \ \ \ \ \vdots  \\
F_n(b) - F_n(a) \\
\end{array}
\right)= 
\left(\begin{array}{ccc}
\frac{\partial F_1}{\partial x_1}(c_1) &\cdots & \frac{\partial F_1}{\partial x_n}(c_1)\\
\  \vdots & & \vdots \\
\frac{\partial F_n}{\partial x_1}(c_n) &\cdots & \frac{\partial F_n}{\partial x_n}(c_n)\\
\end{array}
\right)
\left(\begin{array}{l}
b_1- a_1\\
\ \ \ \ \vdots  \\
b_n -a_n \\
\end{array}
\right).  \]
Since $\det\big(\frac{\partial F_i}{\partial x_j}(c_i)\big)\neq 0$ and $b-a\neq 0$, we conclude that $F(b)\neq F(a)\blacksquare$ 

\section{The Implicit and Inverse Function Theorems.}

\

The first implicit function result we prove regards one equation and several variables. 
We denote the variable in $\mathbb R^{n+1}=\mathbb R^n\times \mathbb R$ by $(x,y)$, where $x=(x_1,\ldots,x_n)$ is in $\mathbb R^n$ and $y$ is in $\mathbb R$. 

\vspace{0,2 cm}
 
\begin{theorem}\label{TEO 1}  Let $F:\Omega \to \mathbb R$ be of class $C^1$ in an open set $\Omega$ inside $\mathbb R^n\times \mathbb R$ and $(a,b)$ be a point in $\Omega$ such that $F(a,b)=0$ and $\frac{\partial F}{\partial y}(a,b)>0$. Then, there exist an open set $X$, inside $\mathbb R^n$ and containing $a$, and an open set $Y$, inside $\mathbb R$ and containing $b$, satisfying the following.
\begin{itemize}
\item[$\bullet$] For each $x$ in $X$ there is a unique $y=f(x)$ in $Y$ such that $F\big(x,f(x)\big)=0$. 
\item[$\bullet$] We have $f(a)=b$. Moreover, $f:X \to Y$ is of class $C^1$ and 
\[\frac{\partial f}{\partial x_j}(x)= - \frac{\frac{\partial F}{\partial x_j}(x,f(x))}{\frac{\partial F}{\partial y}(x,f(x))},\ \textrm{for all}\ x \ \textrm{in}\ X, \ \textrm{where}\ j=1,\ldots,n.\]
\end{itemize}
\end{theorem}
{\bf Proof.} Let us split the proof into three parts.
\begin{itemize}
\item[$\diamond$]{\sf Existence and Uniqueness.} Since $\frac{\partial F}{\partial y}(a,b)>0$, by continuity there exists a non-degenerate $(n+1)$-dimensional parallelepiped $X'\times [b_1,b_2]$, centered at $(a,b)$ and contained in $\Omega$, whose edges are parallel to the coordinate axes such that $\frac{\partial F}{\partial y} > 0$ on $X'\times [b_1,b_2]$. Then, the function $F(a , y)$, where $y$ runs over $[b_1,b_2]$, is strictly increasing and $ F(a ,b )=0$. Thus, we have $F(a ,b_1) < 0$ and $F(a ,b_2 ) > 0$. By the continuity of $F$, there exists an open non-degenerate $n$-dimensional parallelepiped $X$, centered at $a$ and contained in $X'$, whose edges are parallel to the coordinate axes such that for every $x$ in $X$ we have $F(x,b_1 ) <0$ and $F(x, b_2 )>0$. Hence, fixing an arbitrary $x$ in $X$ and employing the intermediate value theorem to the strictly increasing function $ F(x,y)$, where $y$ runs over $[b_1,b_2]$, yields the existence of a unique $y=f(x)$ inside $Y=]b_1,b_2[$ such that $F(x,f(x))=0$. 

\item[$\diamond$]{\sf Continuity.} Let $\overline{b_1 }$ and $\overline{b_2 }$ be such that $ b_1 < \overline{b_1 } < b < \overline{b_2 }< b_2 $. From above, there exists an open set $X''$, contained in $X$ and containing $a$, such that $f(x)$ is in $]\overline{b_1 } ,\overline{b_2 }[$, for all $x$ in $X''$. Thus, $f$ is continuous at $x=a$. Now, given any $a'$ in $X$, we put $b'=f(a')$. Then, $f: X \to Y$ is a solution of the problem $F(x,h(x))=0$, for all $x$ in $X$, with the condition $h(a')= b'$. Thus, from what we have just done it follows that $f$ is continuous at $a'$.

\item[$\diamond$]{\sf Differentiability.} Given $x$ in $X$, let $e_j$ be jth canonical vector in $\mathbb R^n$ and $t \neq 0$ be small enough so that $x+te_j$ belongs to $X$. Putting $P=\big(x,f(x)\big)$ and $Q=\big(x+te_j, f(x+t e_j)\big)$, we notice that $F$ vanishes at $P$ and $Q$. Thus, by employing the mean-value theorem in several variables to $F$ restricted to the segment $\overline{PQ}$ within $X\times]b_1,b_2[$, we find  a point $(\overline{x},\overline{y})$ inside $\overline{PQ}$ and satisfying
\begin{displaymath}
\begin{array}{ll}
0& =F\big(x +te_j, f(x+ te_j)\big) - F\big(x,f(x)\big)=\\
& = \frac{\partial F}{\partial x_j}(\overline{x},\overline{y})t  + \frac{\partial F}{\partial y}(\overline{x},\overline{y})[f(x+te^j) - f(x)].\\
\end{array}
\end{displaymath}
Since $\frac{\partial F}{\partial x_j}$ and $\frac{\partial F}{\partial y}$ are continuous, with $\frac{\partial F}{\partial y}$ not vanishing on $X\times]b_1,b_2[$, the function $f$ is continuous and $(\overline{x},\overline{y})\rightarrow (x,f(x))$ as $t\to 0$. Thus, 
\[\frac{f(x+ te_j) - f(x)}{t} = -\frac{\frac{\partial F}{\partial x_j}(\overline{x},\overline{y})}{\frac{\partial F}{\partial y}(\overline{x},\overline{y})}\rightarrow - \frac{\frac{\partial F}{\partial x_j}(x,f(x))}{\frac{\partial F}{\partial y}(x,f(x))}\ \textrm{as}\ t \to 0.\]
This gives the desired formula for $\frac{\partial f}{\partial x_j}$ and implies that $f$ is of class $C^1\blacksquare$
\end{itemize}

\vspace{0,2 cm}

Next, we prove the general implicit function theorem. In general, we apply this theorem when we have a nonlinear system with $m$ equations and $n+m$ variables. Analogously to a linear system, we interpret $n$ variables as independent variables  and determine the remaining $m$ variables, called dependent variables, as a function of the $n$ independent variables. 

Let us introduce some helpful notation. As before, we denote by $x=(x_1,...,x_n)$ a point in $\mathbb R^n$ and by $y=(y_1,\ldots,y_m)$ a point in $\mathbb R^m$. 
Given $\Omega$ an open subset of $\mathbb R^n\times \mathbb R^m$ and a  differentiable function $F:\Omega \to \mathbb R^m$, we write   $F=(F_1,\ldots,F_m)$, with $F_i$ the ith component of $F$ and $i=1,\ldots, m$. We put
\[\frac{\partial F}{\partial y}=\left(\frac{\partial F_i}{\partial y_j}\right)_{\substack{1\leq i\leq m\\1 \leq j\leq m}}=
\left(\begin{array}{lllll}
\frac{\partial F_1}{\partial y_1} & \cdots & \frac{\partial F_1}{\partial y_m}\\
\ \ \ \vdots  &   & \ \ \ \vdots \ \ \      \\
\frac{\partial F_m}{\partial y_1} & \cdots & \frac{\partial F_m}{\partial y_m}
\end{array}
\right).\]
Analogously, we define the matrix $\frac{\partial F}{\partial x}=\big(\frac{\partial F_i}{\partial x_k}\big)$, where $1\leq i \leq m$ and $1\leq k \leq n$.

\begin{theorem}\label{TEO 2}{\bf (The Implicit Function Theorem).}  Let $F$ be in $C^1(\Omega;\mathbb R^m)$, with $\Omega$ an open set in $\mathbb R^n\times \mathbb R^m$, and $(a,b)$ a point in $\Omega$ such that $F(a,b)=0$ and $\frac{\partial F}{\partial y}(a,b)$ is invertible. Then, there exist an open set $X$, inside $\mathbb R^n$ and containing $a$, and an open set $Y$, inside $\mathbb R^m$ and containing $b$, satisfying the following.
\begin{itemize}
\item[$\bullet$] For each $x$ in $X$, there is a unique $y=f(x)$ in $Y$ such that $F\big(x,f(x)\big)=0$. 
\item[$\bullet$] We have $f(a)=b$. Moreover, $f:X \to Y$ is of class $C^1$ and 
\[Jf(x)  = - \left[\frac{\partial F}{\partial y }(x,f(x))\right]_{m\times m}^{-1}\left[\frac{\partial F}{\partial x}(x,f(x))\right]_{m\times n},\ \textrm{for all}\ x \ \textrm{in}\ X.\]
\end{itemize}
\end{theorem}
{\bf Proof.} We split the proof into four parts.
\begin{itemize}
\item[$\diamond$] {\sf Finding $Y$.} Defining $\Phi(x,y)=\big(x,F(x,y))$, where $(x,y)$ is in $\Omega$, we have
$$\det J\Phi=\det\left(\begin{array}{l|l}
\ I  &\ 0\\ 
\hline
\frac{\partial F}{\partial x}  & \frac{\partial F}{\partial y}
\end{array} \right)= \det \frac{\partial F}{\partial y},$$ 
with $I$ the identity matrix of order $n$ and $0$ the $n\times m$ zero matrix. Thus, shrinking $\Omega$ if needed, by Lemma 4 we may and do assume, without any loss of generality, that $\Phi$ is injective and  $\Omega=X'\times Y$, with $X'$ an open set in $\mathbb R^n$ that contains $a$ and $Y$ an open set in $\mathbb R^m$ that contains $b$.

\item[$\diamond$]{\sf Existence and differentiability.} We claim that the equation $F(x,f(x))=0$, with the condition $f(a)=b$, has a solution $f=f(x)$ of class $C^1$ on some open set containing $a$. Let us prove it by induction on $m$. 

The case $m=1$ follows from Theorem 1. Let us assume that the claim is true for $m-1$. Then, for the case $m$, following the notation $F=(F_1,\ldots,F_m)$ we write $\mathcal{F}=(F_2,\ldots,F_m)$.  Furthermore, we put $(x;y)=(x_1,\ldots x_n;y_1,\ldots, y_m)$, $y'=(y_2,\ldots,y_m)$, $y=(y_1;y')$, and $(x;y)=(x;y_1;y')$.

Let us consider the invertible matrix $J=\frac{\partial F}{\partial y}(a,b)$ and the associated bijective linear function $\mathcal{J}:\mathbb R^m \to \mathbb R^m$.
By Lemma 2, the function $G(x;z) = F[x;b + \mathcal{J}^{-1}(z-b)]$, defined in some open subset of $\mathbb R^n\times \mathbb R^m$ that contains $(a,b)$,  satisfy $\frac{\partial G}{\partial z}(a;b)=JJ^{-1}$ and the condition $G(a;b)=0$. Hence, we may assume that $J$ is the identity matrix of order $m$.

Now, let us consider the equation $F_1(x;y_1;y')=0$ with the condition $y_1(a;b')=b_1$. Since $\frac{\partial F_1}{\partial y_1}(a;b_1;b')=1$, there exists a function $y_1=\varphi(x;y')$ of class $C^1$ on some open set containing $(a;b')$ that satisfies 
\[F_1[x;\varphi(x;y');y']=0 \  \textrm{and the condition}\ \varphi(a;b')=b_1.\]
Next, substituting $y_1=\varphi(x;y')$ into  $\mathcal{F}(x;y_1;y')=0$, we look at solving
\[\textrm{the equation} \ \mathcal{F}[x;\varphi(x,y');y']=0,\ \textrm{with the condition}\ y'(a)=b'.\]
Differentiating $\mathcal{F}[x;\varphi(x;y');y']$, with respect to $y_2,\ldots,y_m$, we find
\[\frac{\partial F_i}{\partial y_1}(a;b)\frac{\partial \varphi}{\partial y_j}(a;b') + \frac{\partial F_i}{\partial y_j}(a;b)=\ 0 + \frac{\partial F_i}{\partial y_j}(a;b), \  \textrm{where}\ 2\leq i,j\leq m.\] 
The matrix $\left(\frac{\partial F_i}{\partial y_j}(a;b)\right)$, where $2\leq i,j\leq m$, is the identity one of order $m-1$. Hence, by induction hypothesis there exists a function $\psi$ of class $C^1$ on an open set $X$ containing $a$ that satisfies 
\[ \mathcal{F}\big[x;\varphi\big(x;\psi(x)\big),\psi(x)\big]=0, \ \textrm{for all}\ x \ \textrm{in}\ X, \ 
\textrm{and the condition} \ \psi(a)=b'.\] 
We also have $F_1\big[x;\varphi\big(x;\psi(x)\big);\psi(x)\big]=0$, for all $x$ in $X$. Defining $f(x) = \big(\varphi(x;\psi(x));\psi(x)\big)$, with $x$ in $X$, we obtain $F[x;f(x)]=0$, for all $x$ in $X$, and $f(a)=\big(\varphi(a;b');b'\big)=(b_1;b')=b$, where $f$ is of class $C^1$ on $X$.
\item[$\diamond$] {\sf Differentiation formula.} Differentiating $F[x,f(x)]=0$ we find
\[\frac{\partial F_i}{\partial x_k} + \sum_{j=1}^m\frac{\partial F_i}{\partial y_j}\frac{\partial f_j}{\partial x_k}=0,\ \textrm{with}\ 1\leq i\leq m\ \textrm{and}\ 1\leq k\leq n.\]
In matricial form, we write $\frac{\partial F}{\partial x}\big(x,f(x)\big) + \frac{\partial F}{\partial y}\big(x,f(x)\big)Jf(x)=0$.

\item[$\diamond$] {\sf Uniqueness.} Let $g$ satisfy $F\big(x,g(x)\big)=0$, for all $x$ in $X$, and $g(a)=b$. Thus, we have $\Phi(x,g(x)\big)= \big(x,F(x,g(x)\big)=(x,0)= \Phi\big(x,f(x)\big)$, for all $x$ in $X$. The injectivity of $\Phi$ implies that $g(x)=f(x)$, for all $x$ in $X$$\blacksquare$
\end{itemize}

\begin{theorem} \label{TEO 3}{\bf (The Inverse Function Theorem).} Let $F:\Omega \to \mathbb R^n$, where $\Omega$ is an open set in $\mathbb R^n$, be of class $C^1$ and $p$ a point in $\Omega$ such that $JF(p)$ is invertible. Then, there exist an open set $X$ containing $p$, an open set $Y$ containing $F(p)$, and a function  $G:Y\to X$ of class $C^1$ that satisfies $F\big(G(y)\big)=y$, for all $y$ in $Y$, and $G\big(F(x)\big)=x$, for all $x$ in $X$. Moreover, 
\[JG(y)= JF\big(G(y)\big)^{-1}, \ \textrm{for all} \ y \ \textrm{in}\ Y.\]
\end{theorem}
{\bf Proof.} We split the proof into two parts: existence and differentiation formula.

\begin{itemize}
\item[$\diamond$] {\sf Existence.} Shrinking $\Omega$, if necessary, by Lemma 4 we may assume that $F$ is injective. The function $\Phi(x,y)= F(x)-y$, where $(x,y)$ is in $\Omega\times \mathbb R^n$, is of class $C^1$ and satisfies $\frac{\partial \Phi}{\partial x}\big(F(p),p\big)=JF(p)$. From the implicit function theorem it follows that there exists an open set $Y$ containing $F(p)$ and a  function $G:Y \to \Omega$ of class $C^1$ such that $\Phi\big(G(y),y\big) =F\big(G(y)\big) - y=0$, for all $y$ in $Y$. That is, we have $F\big(G(y)\big)=y$, for all $y$ in $Y$. 

Hence, the set $Y$ is contained in the image of $F$. Since $F$ is continuous and injective,  the pre-image set $X=F^{-1}(Y)$ is open, contains $p$, and $F$ maps $X$ bijectively onto $Y$.

The identity $F\big(G(y)\big)=y$, for all $y$ in $Y$, implies that $G$ maps $Y$ to $X$. Since $F$ is bijective from $X$ to $Y$, the function $G$ is bijective from $Y$ to $X$.

\item[$\diamond$] {\sf Differentiation formula.} Let us write $F(x)=\big(F_1(x),\ldots,F_n(x)\big)$ and $G(y)=\big(G_1(y),\ldots,G_n(y)\big)$. Differentiating  $\big(G_1(F(x)),\ldots,G_n(F(x))\big)$ we obtain
\[ \sum_{k=1}^n\frac{\partial G_i}{\partial y_k}\frac{\partial F_k}{\partial x_j}= \frac{\partial x_i}{\partial x_j}=
\left\{\begin{array}{ll}
1,\ \textrm{if} \ i=j,\\
0, \ \textrm{if}\ i\neq j\,\blacksquare
\end{array}
\right.\]
\end{itemize}

\paragraph{Acknowledgments.}  I thank Professors Robert B. Burckel, Paulo A. Martin, and A. Lymberopoulos for their very valuable comments and suggestions.

\

\bigskip

\bigskip

\noindent\textit{Departamento de Matemática,
Universidade de São Paulo\\
Rua do Matão 1010 - CEP 05508-090\\
São Paulo, SP - Brasil\\
oliveira@ime.usp.br}

\bigskip

\end{document}